\documentclass[proceedings]{aofa}
\pdfoutput = 1
\newif\ifprivate

\ifprivate
\usepackage[mark]{gitinfo2}
\renewcommand{\gitMark}{\jobname\,\textbullet{}\,\gitFirstTagDescribe\,\textbullet{}\,\gitAuthorName,\,\gitAuthorIsoDate}
\newcommand{\TODO}[1]%
{\par\fbox{\begin{minipage}{0.9\linewidth}\textbf{TODO:} #1\end{minipage}}\par}
\fi

\usepackage{amsthm}

\usepackage{amsmath}
\usepackage{amssymb}
\usepackage{hyperref}
\usepackage[utf8]{inputenc}
\usepackage{mathtools}

\newcommand{\bfr}{\mathbf{r}}

\newcommand{\bfOmega}{\boldsymbol{\Omega}}
\newcommand{\bfs}{\mathbf{s}}

\newcommand{\bfQ}{\mathbf{Q}}
\newcommand{\bfk}{\mathbf{k}}
\newcommand{\bft}{\mathbf{t}}

\newcommand{\bfx}{\mathbf{x}}
\newcommand{\bfX}{\mathbf{X}}
\newcommand{\bfy}{\mathbf{y}}
\newcommand{\bfY}{\mathbf{Y}}
\newcommand{\bfz}{\mathbf{z}}
\newcommand{\bfZ}{\mathbf{Z}}
\newcommand{\bfzero}{\boldsymbol{0}}
\newcommand{\bfones}{\boldsymbol{1}}

\newcommand{\C}{\mathbb{C}}

\newcommand{\E}{\mathbb{E}}
\DeclareMathOperator{\grad}{grad}

\newcommand{\MR}[1]{}
\renewcommand{\P}{\mathbb{P}}
\newcommand{\R}{\mathbb{R}}

\newcommand{\sphi}{\sqrt{\phi_{n}}}
\newcommand{\sphii}{\phi_n^{-1/2}}
\newcommand{\stirlingpartition}[2]{\genfrac{\{}{\}}{0pt}{}{#1}{#2}}

\DeclareMathOperator{\Cov}{Cov}
\newcommand{\cupdot}{\mathbin{\mathaccent\cdot\cup}}

\DeclarePairedDelimiter{\abs}{\lvert}{\rvert}
\DeclarePairedDelimiter{\norm}{\lVert}{\rVert}
\DeclarePairedDelimiter{\iverson}{[}{]}

\newtheorem{theorem}{Theorem}
\newtheorem*{theoremnonumber}{Theorem}
\newtheorem{lemma}{Lemma}[section]
\newtheorem{proposition}[lemma]{Proposition}
\newtheorem{corollary}[lemma]{Corollary}

\theoremstyle{definition}
\newtheorem{definition}[lemma]{Definition}

\theoremstyle{remark}

\newenvironment{sketch}[1][]{}{}
\def\sketch[#1]{\noindent\textbf{#1:} \ignorespaces}

\numberwithin{equation}{section}

\author[Clemens Heuberger and Sara Kropf]{Clemens Heuberger\addressmark{1}\thanks{The authors are supported by the Austrian Science Fund (FWF):
  P~24644-N26.} and Sara Kropf\addressmark{1,2}\addressmark{\dag}}

\title{On the Higher Dimensional Quasi-Power Theorem and a
  Berry--Esseen Inequality}

\address{\addressmark{1}Institut f\"ur Mathematik, Alpen-Adria-Universit\"at
  Klagenfurt, Austria, \{clemens.heuberger, sara.kropf\}@aau.at\\
\addressmark{2}Institute of Statistical Science, Academia
  Sinica, Taipei, Taiwan, sarakropf@stat.sinica.edu.tw}

\begin{document}
\maketitle
\begin{abstract}Hwang's quasi-power theorem asserts that a sequence of
  random variables whose moment generating functions are approximately given by
  powers of some analytic function is asymptotically normally distributed.
  This theorem is generalised to higher dimensional random
  variables. To obtain this result, a higher dimensional analogue of the Berry--Esseen
  inequality is proved, generalising a two-dimensional version by
  Sadikova.
\end{abstract}

\section{Introduction}
Asymptotic normality is a frequently occurring phenomenon in combinatorics, the
classical central limit theorem being the very first example. The first step in
the proof is the observation that the
moment generating function of the sum of $n$ identically independently distributed
random variables is the $n$-th power of the moment generating function of the
distribution underlying the summands. As similar moment generating functions
occur in many examples in combinatorics, a general theorem to prove asymptotic
normality is desirable. Such a theorem was proved by
Hwang~\cite{Hwang:1998}, usually called the ``quasi-power theorem''.

\begin{theoremnonumber}[Hwang~\cite{Hwang:1998}]
  Let $\{\Omega_n\}_{n\ge 1}$ be a sequence of  integral random
  variables. Suppose that the moment generating function satisfies the asymptotic expression
  \begin{equation}\label{eq:moments-1d}
    M_n(s):=\E(e^{\Omega_ns})=e^{W_n(s)}(1+O(\kappa_n^{-1})),
  \end{equation}
  the $O$-term being uniform for $\abs{s}\le \tau$, $s\in\C$,
  $\tau>0$, where
  \begin{enumerate}
  \item $W_n(s)=u(s)\phi_{n}+v(s)$, with $u(s)$ and $v(s)$
    analytic for $\abs{s}\le \tau$ and independent of $n$; and $u''(0)\neq 0$;
  \item $\lim_{n\to\infty}\phi_{n}=\infty$;
  \item $\lim_{n\to\infty}\kappa_n=\infty$.
  \end{enumerate}

  Then the distribution of $\Omega_n$ is asymptotically normal, i.e.,
  \begin{equation*}
    \sup_{x\in\R}\left\vert\P\left(\frac{\Omega_n- u'(0)\phi_{n}}{\sqrt{u''(0)\phi_{n}}} <
      x\right)-
    \Phi(x)\right\vert=O\left(\frac{1}{\sqrt{\phi_{n}}}+\frac{1}{\kappa_n}\right),
  \end{equation*}
  where $\Phi$ denotes the standard normal distribution
  \begin{equation*}
    \Phi(x)=\frac{1}{\sqrt{2\pi}}\int_{-\infty}^{x}\exp\left(-\frac12
      y^2\right)\,dy.
  \end{equation*}
\end{theoremnonumber}

See Hwang's article~\cite{Hwang:1998} as well as
Flajolet-Sedgewick~\cite[Sec.~IX.5]{Flajolet-Sedgewick:ta:analy} for many
applications of this theorem.
A generalisation of the quasi-power theorem to dimension~$2$ has been provided
 in \cite{Heuberger:2007:quasi-power}. It has been used in \cite{Heuberger-Prodinger:2006:analy-alter},
 \cite{Heuberger-Prodinger:2007:hammin-weigh},
 \cite{Eagle-Gao-Omar-Panario:2008:distr:short},
 \cite{Heuberger-Kropf-Wagner:2014:combin-charac} and
 \cite{Kropf:2015:varian-and}. In~\cite[Thm.~2.22]{Drmota:2009:random},
 an $m$-dimensional version of the quasi-power theorem is stated
 without speed of convergence. Also
 in~\cite{Bender-Richmond:1983:centr}, such an $m$-dimensional theorem
 without speed of convergence is proved. There, several multidimensional applications
 are given, too. 

In contrast to many results about the speed of convergence in classical probability
 theory (see, e.g.,~\cite{Gut:2005:probab}), the sequence of
 random variables is not assumed to be independent. The only assumption is that the
 moment generating function behaves asymptotically like a large
 power. This mirrors the fact that the moment generating function of
 the sum of independent, identically distributed random variables is exactly a
 large power. The advantage is that the asymptotic expression~\eqref{eq:moments-1d} arises naturally in combinatorics by
 using techniques such as singularity analysis or saddle point
 approximation (see \cite{Flajolet-Sedgewick:ta:analy}).  

The purpose of this article is to generalise the quasi-power theorem
including the speed of convergence to
arbitrary dimension $m$. We first state this main result in
Theorem~\ref{th:quasi-power-dD} in this section. In
Section~\ref{sec:Berry--Esseen}, a new Berry--Esseen inequality
(Theorem~\ref{theorem:Berry-Esseen-dimension-m}) is presented, which we
use to prove the $m$-dimensional quasi-power theorem. We give sketches of the
proofs of these two theorems in
Section~\ref{sec:sketch-proofs}. All details of these proofs can be
found in the full version of this extended abstract. In
Section~\ref{sec:exampl-mult-centr}, we give some applications of the
multidimensional quasi-power theorem.

We use the following conventions: vectors are denoted by
boldface letters such as $\bfs$, their components are then denoted by
regular letters with indices such as $s_j$. For a vector $\bfs$, $\|\bfs\|$ denotes the
maximum norm $\max\{\abs{s_j}\}$.
All implicit constants of
$O$-terms may depend on the dimension $m$ as well as on $\tau$ which is introduced in
Theorem~\ref{th:quasi-power-dD}.


Our first main result is the following $m$-dimensional version of Hwang's theorem.

\begin{theorem}\label{th:quasi-power-dD}
  Let $\{\bfOmega_n\}_{n\ge 1}$ be a sequence of $m$-dimensional real random
  vectors. Suppose that the moment generating function satisfies the asymptotic expression
  \begin{equation}\label{eq:moment-asymp}
    M_n(\bfs):=\E(e^{\langle \bfOmega_n,\bfs\rangle})=e^{W_n(\bfs)}(1+O(\kappa_n^{-1})),
  \end{equation}
  the $O$-term being uniform for $\norm{\bfs}\le \tau$, $\bfs\in\C^m$,
  $\tau>0$, where
  \begin{enumerate}
  \item $W_n(\bfs)=u(\bfs)\phi_{n}+v(\bfs)$, with $u(\bfs)$ and $v(\bfs)$
    analytic for $\norm{\bfs}\le \tau$ and independent of $n$; and the Hessian
    $H_u(\bfzero)$ of $u$ at the origin is non-singular;
  \item $\lim_{n\to\infty}\phi_{n}=\infty$;
  \item $\lim_{n\to\infty}\kappa_n=\infty$.
  \end{enumerate}

  Then, the distribution of $\bfOmega_n$ is
  asymptotically normal with speed of convergence $O(\sphii)$, i.e.,
  \begin{equation}\label{eq:quasi-power-result}
    \sup_{\bfx\in\R^{m}}\left\vert\P\left(\frac{\bfOmega_n-\grad u (\bfzero)\phi_{n}}{\sqrt{\phi_{n}}} \le
      \bfx\right)-
    \Phi_{H_u(\bfzero)}(\bfx)\right\vert=O\left(\frac{1}{\sqrt{\phi_{n}}}\right),
  \end{equation}
  where $\Phi_{\Sigma}$ denotes the distribution function
  of the non-degenerate $m$-dimensional normal distribution
  with mean $\bfzero$ and variance-covariance
  matrix $\Sigma$, i.e.,
  \begin{equation*}
    \Phi_\Sigma(\bfx)=\frac{1}{(2\pi)^{m/2}\sqrt{\det \Sigma}}\int_{\bfy\le \bfx}\exp\left(-\frac12
      \bfy^\top \Sigma^{-1} \bfy \right)\,d\bfy,
  \end{equation*}
  where $\bfy\le \bfx$ means $y_\ell\le x_\ell$ for $1\le \ell\le m$.

  If $H_{u}(\bfzero)$ is singular, the random variables
  \begin{equation*}
    \frac{\bfOmega_{n}-\grad u(\bfzero)\phi_{n}}{\sqrt{\phi_{n}}}
  \end{equation*}
  converge in distribution to a
  degenerate normal distribution with mean $\bfzero$ and variance-covariance
  matrix $H_{u}(\bfzero)$.
\end{theorem}

Note that in the case of the singular $H_{u}(\bfzero)$, a uniform speed of
convergence cannot be guaranteed.
 To see this,
  consider the (constant) sequence of random variables $\Omega_{n}$
  which takes values $\pm1$ each with probability $1/2$. Then the
  moment generating function is $(e^{t}+e^{-t})/2$,
which is of the form \eqref{eq:moment-asymp} with $\phi_{n}=n$,
$u(s)=0$, $v(s)=\log (e^{t}+e^{-t})/2$ and $\kappa_{n}$ arbitrary. However, the distribution function of $\Omega_{n}/\sqrt{n}$ is
  given by
  \begin{equation*}
    \mathbb{P}\biggl(\frac{\Omega_{n}}{\sqrt{n}}\le x\biggr)=
    \begin{cases}
      0& \text{if }x<-1/\sqrt{n},\\
      1/2& \text{if }-1/\sqrt{n}\le x<1/\sqrt{n},\\
      1& \text{if }1/\sqrt{n}\le x,
    \end{cases}
  \end{equation*}
  which does not converge uniformly.

In contrast to the original quasi-power theorem, the error term in our result does not
contain the summand $O(1/\kappa_n)$. In fact, this summand could also be
omitted in the original proof of the quasi-power theorem by using a better
estimate for the error $E_{n}(\bfs)=M_n(\bfs)e^{-W_{n}(\bfs)}-1$.

The proof of Theorem~\ref{th:quasi-power-dD} relies on an
$m$-dimensional Berry--Esseen inequality
(Theorem~\ref{theorem:Berry-Esseen-dimension-m}). It is a generalisation of
Sadikova's result~\cite{Sadikova:1966:esseen, Sadikova:1966:esseen:englisch} in
dimension $2$. The main challenge is to provide a version which leads to
bounded integrands around the origin, but still allows to use excellent bounds for the
tails of the characteristic functions. To achieve this, linear combinations
involving all partitions of the set $\{1,\ldots, m\}$ are used.

Note that there are several generalisations of the one-dimensional
Berry--Esseen inequality \cite{Berry:1941:gauss,Esseen:1945:fourier}
to arbitrary dimension, see, e.g., Gamkrelidze~\cite{Gamkrelidze:1977,
  Gamkrelidze:1977:englisch} and Prakasa
Rao~\cite{Rao:2002:anoth-esseen}. However, using these results would lead to
the less precise error term in~\eqref{eq:quasi-power-result}, see the end of
Section~\ref{sec:Berry--Esseen} for more details. For that reason we generalise Sadikova's
result, which was already successfully used by the first author
in~\cite{Heuberger:2007:quasi-power} to prove a $2$-dimensional
quasi-power theorem. Also note
 that our theorem can deal with discrete random variables, in contrast
 to \cite{Roussas:2001:esseen}, where density functions are considered.

For the sake of completeness, we also state the following result
about the moments of $\bfOmega_{n}$.
\begin{proposition}\label{proposition:moments}
  The cross-moments of $\bfOmega_{n}$ satisfy
  \begin{equation*}
    \frac{1}{\prod_{\ell=1}^{m}k_{\ell}!}\mathbb E\Big(\prod_{\ell=1}^{m}\Omega_{n,\ell}^{k_{\ell}}\Big)=p_{\bfk}(\phi_{n})+O\Big(\kappa_{n}^{-1}\phi_{n}^{k_{1}+\cdots+k_{m}}\Big),
  \end{equation*}
for $k_{\ell}$ nonnegative integers, where $p_{\bfk}$ is a polynomial of degree
$\sum_{\ell=1}^{m}k_{\ell}$ defined by
\begin{equation*}
  p_{\bfk}(X)=[s_{1}^{k_{1}}\cdots s_{m}^{k_{m}}]e^{u(\bfs)X+v(\bfs)}.
\end{equation*}

In particular, the mean and the variance-covariance matrix are
\begin{align*}
  \mathbb E(\bfOmega_{n})&=\grad u(\bfzero)\phi_{n}+\grad
  v(\bfzero)+O(\kappa_{n}^{-1}),\\
  \Cov(\bfOmega_{n})&=H_{u}(\bfzero)\phi_{n}+H_{v}(\bfzero)+O(\kappa_{n}^{-1}),
\end{align*}
respectively.
\end{proposition}

\section{A Berry--Esseen Inequality}\label{sec:Berry--Esseen}

This section is devoted to a generalisation of Sadikova's Berry--Esseen
inequality~\cite{Sadikova:1966:esseen, Sadikova:1966:esseen:englisch} in
dimension 2 to dimension $m$. Before stating the theorem, we introduce our notation.

Let $L=\{1,\ldots, m\}$. For $K\subseteq L$, we write
$\bfs_K=(s_k)_{k\in K}$ for the projection of $\bfs\in\C^L$ to $\C^K$.
For $J\subseteq K\subseteq L$, let
$\chi_{J,K}\colon \C^{J}\to\C^{K}$, $(s_{j})_{j\in J}\mapsto
(s_{k}\iverson{k\in J})_{k\in K}$ be an injection from $\C^{J}$ into
$\C^{K}$. Similarly, let $\psi_{J,K}\colon \C^{K}\to\C^{K}$,
$(s_{k})_{k\in K}\mapsto (s_{k}\iverson{k\in J})_{k\in K}$ be the
projection which sets all coordinates corresponding to $K\setminus J$ to
$0$.

We denote the
  set of all partitions of $K$ by $\Pi_K$. We consider a partition as a set
  $\alpha=\{J_{1},\ldots,J_{k}\}$. Thus $\abs{\alpha}$ denotes the number of parts of the partition
$\alpha$. Furthermore, $J\in\alpha$ means that $J$ is a part of the
partition $\alpha$.

Now, we can define an operator which we later use to state our Berry--Esseen
inequality. The motivation behind this definition is explained at the end of
this section.

\begin{definition}\label{definition:Lambda-K}
  Let $K\subseteq L$ and $h\colon \C^K\to \C$. 
 We define the non-linear operator
  \begin{equation*}
    \Lambda_K(h):=\sum_{\alpha\in\Pi_K}\mu_\alpha \prod_{J\in
      \alpha}h\circ \psi_{J, K}
  \end{equation*}
  where
  \begin{equation*}
    \mu_\alpha = (-1)^{\abs{\alpha}-1}(\abs{\alpha}-1)!\,.
  \end{equation*}
We denote $\Lambda_{L}$ briefly by $\Lambda$.
\end{definition}

For any random variable $\bfZ$, we denote its cumulative distribution function
by $F_\bfZ$ and its
characteristic function by $\varphi_\bfZ$.

With these definitions, we are able to state our second main result, an
$m$-dimensional version of the Berry--Esseen inequality.

\begin{theorem}\label{theorem:Berry-Esseen-dimension-m}
Let $m\ge 1$ and $\bfX$ and $\bfY$ be $m$-dimensional random variables. Assume
that $F_\bfY$ is differentiable.

Let
\begin{align*}
  A_j&=\sup_{\bfy \in\R^m}\frac{\partial F_\bfY(\bfy)}{\partial y_j},\\
  B_j&=\sum_{k=1}^{j} \stirlingpartition{j}{k} k!\ ,\\
  C_1&=\sqrt[3]{\frac{32}{\pi\bigl(1-\bigl(\frac{3}{4}\bigr)^{1/m}\bigr)}},\\
  C_2&=\frac{12}{\pi}
\end{align*}
for $1\le
  j\le m$ where $\stirlingpartition{j}{k}$ denotes a Stirling partition number
(Stirling number of the second kind).

Let $T>0$ be fixed. Then
\begin{equation}\label{eq:Berry-Esseen}
  \begin{aligned}
    \sup_{\bfz\in\R^m}\abs{F_{\bfX}(\bfz)-F_{\bfY}(\bfz)}&\le
    \frac{2}{(2\pi)^m} \int_{\norm{\bft}\le
      T}\abs[\Big]{\frac{\Lambda(\varphi_{\bfX})(\bft)-\Lambda(\varphi_{\bfY})(\bft)}{\prod_{\ell\in
          L} t_\ell}}\,d\bft \\
    &\qquad+ 2\sum_{\emptyset\neq J\subsetneq
      L}B_{m-\abs{J}}\sup_{\bfz_J\in\R^J}\abs[\big]{F_{\bfX_{J}}(\bfz_J)-F_{\bfY_{J}}(\bfz_J)}
    \\
    &\qquad +\frac{2\sum_{j=1}^m A_j}{T}(C_1+C_2).
  \end{aligned}
\end{equation}
Existence of $\E(\bfX)$ and $\E(\bfY)$ is sufficient for the finiteness
of the integral in \eqref{eq:Berry-Esseen}.
\end{theorem}
Let us give two remarks on the distribution functions occurring in this theorem: The distribution function $F_\bfY$ is non-decreasing in every
variable, thus $A_j>0$ for all $j$. Furthermore, our general notations imply that $F_{\bfX_J}$ is a marginal
distribution of $\bfX$.

The numbers $B_j$ are known as ``Fubini numbers'' or ``ordered Bell numbers''.
They form the sequence \href{http://oeis.org/A000670}{A000670} in \cite{OEIS:2015}.

Recursive application of \eqref{eq:Berry-Esseen} leads to the following
corollary, where we no longer explicitly state the constants depending on the
dimension.

\begin{corollary}\label{corollary:Berry-Esseen}
  Let $m\ge 1$ and $\bfX$ and $\bfY$ be $m$-dimensional random
  variables. Assume
that $F_\bfY$ is differentiable and
  let
  \begin{equation*}
    A_j=\sup_{\bfy \in\R^m}\frac{\partial F_\bfY(\bfy)}{\partial y_j}, \qquad 1\le
    j\le m.
  \end{equation*}

  Then
  \begin{multline}\label{eq:Berry-Esseen-recursive}
      \sup_{\bfz\in\R^m}\abs{F_{\bfX}(\bfz)-F_{\bfY}(\bfz)}\\=
      O\biggl(\sum_{\emptyset \neq K\subseteq L}\int_{\norm{\bft_K}\le
        T}\abs[\Big]{\frac{\Lambda_K(\varphi_{\bfX}\circ\chi_{K, L})(\bft_K)-\Lambda_K(\varphi_{\bfY}\circ\chi_{K, L})(\bft_K)}{\prod_{k\in
            K} t_k}}\,d\bft_K + \frac{\sum_{j=1}^m A_j}{T}\biggr)
  \end{multline}
  where the $O$-constants only depend on the dimension $m$.

  Existence of $\E(\bfX)$ and $\E(\bfY)$ is sufficient for the finiteness
  of the integrals in \eqref{eq:Berry-Esseen-recursive}.
\end{corollary}

In order to explain the choice of the operator $\Lambda$, we first state it in
dimension $2$:
\begin{equation}\label{eq:Sadikova-simple}
  \Lambda(h)(s_1, s_2) = h(s_1, s_2) - h(s_1, 0)h(0, s_2).
\end{equation}
This coincides with Sadikova's definition. This also shows that our
operator is non-linear as, e.g., $\Lambda(s_{1}+s_{2})(s_{1},s_{2})\neq\Lambda(s_{1})(s_{1},s_{2})+\Lambda(s_{2})(s_{1},s_{2})$.

In Theorem~\ref{theorem:Berry-Esseen-dimension-m}, we apply $\Lambda$ to
characteristic functions; so we may restrict our attention to functions $h$
with $h(\bfzero)=1$. From~\eqref{eq:Sadikova-simple}, we see that
$\Lambda(h)(s_1, 0) = \Lambda(h)(0, s_2)=0$, so that $\Lambda(h)(s_1,
s_2)/(s_1s_2)$ is bounded around the origin. This is essential for the boundedness
of the integral in Theorem~\ref{theorem:Berry-Esseen-dimension-m}. In general,
this property will be guaranteed by our particular choice of coefficients. It
is no coincidence that for $\alpha\in \Pi_L$, the coefficient $\mu_\alpha$
equals the value  $\mu(\alpha, \{L\})$ of the Möbius function in the lattice of partitions:
Weisner's theorem (see Stanley~\cite[Corollary~3.9.3]{Stanley:2012:enumer_1}) is crucial in the proof that $\Lambda(h)(\bfs)/(s_1\cdots
s_m)$ is bounded around the origin.

The second property is that our proof of the quasi-power theorem needs
estimates for the tails of the integral in Theorem~\ref{theorem:Berry-Esseen-dimension-m}. These estimates have to be
exponentially small in every variable, which means that every variable has to
occur in every summand. This is trivially fulfilled as every summand
in the definition of $\Lambda$ is
formulated in terms of a partition.

Note that Gamkrelidze~\cite{Gamkrelidze:1977:englisch} (and also
Prakasa Rao~\cite{Rao:2002:anoth-esseen}) use a linear operator $L$ mapping $h$ to
\begin{equation}\label{eq:gamkrelidze}
 (s_1, s_2) \mapsto h(s_1, s_2) - h(s_1, 0) - h(0, s_2).
\end{equation}
When taking the difference of two characteristic functions, we may assume that
$h(0, 0)=0$ so that the first crucial property as defined above still
holds. However, the tails are no longer exponentially small in every
variable: The last summand $h(0,s_{2})$ in \eqref{eq:gamkrelidze} is not exponentially small in
$s_{1}$ because it is independent of $s_{1}$ and nonzero in
general. However, the first two summands are exponentially small in
$s_{1}$ by our assumption~\eqref{eq:moment-asymp}.

For that reason, using the Berry--Esseen inequality by
Gamkrelidze~\cite{Gamkrelidze:1977:englisch} to prove a quasi-power theorem leads to a less precise error term 
$O(\phi_{n}^{-1/2}\log^{m-1}\phi_n)$
in~\eqref{eq:quasi-power-result}. It can be shown that the less precise
error term necessarily appears when using Gamkrelidze's result by considering
the example of $\bfOmega_n$ being the $2$-dimensional vector consisting of a normal distribution
with mean $-1$ and variance $n$ and a normal distribution with mean $0$ and
variance $n$. This is a consequence of the linearity of the
operator $L$ in Gamkrelidze's result.

\section{Examples of Multidimensional Central Limit
  Theorems}\label{sec:exampl-mult-centr}
In this section, we give two examples from combinatorics where we can
apply Theorem~\ref{th:quasi-power-dD}. Asymptotic normality was
already shown in earlier publications
\cite{Drmota:1997:system-funct-equat,Bender-Richmond:1983:centr}, but
we additionally provide an estimate for the
speed of convergence.
\subsection{Context-Free Languages}
Consider the following example of a context-free grammar $G$ with
non-terminal symbols $S$ and $T$, terminal symbols $\{a,b,c\}$,
starting symbol $S$ and the rules 
\begin{equation*}
  P=\{S\to aSbS,\, S\to bT,\, T\to bS,\, T\to cT,\, T\to a\}.
\end{equation*}
The corresponding context-free language $L(G)$ consists of all words which
can be generated starting with $S$ using the rules in $P$ to replace
all non-terminal symbols. For example, $abcabababba\in L(G)$ because it can be derived as
\begin{equation*}
  S\to aSbS
  \to abTbaSbS
\to abcTbabTbbT
\to abcabababba.
\end{equation*}

Let $\P(\bfOmega_{n}=\bfx)$ be the probability that a word of length $n$
in $L(G)$ consists of $x_{1}$ and $x_{2}$ terminal symbols
$a$ and $b$, respectively. Thus there are $n-x_{1}-x_{2}$ terminal
symbols $c$. For simplicity, this random variable is only
$2$-dimensional. But it can be easily
extended to higher dimensions.

Following
Drmota~\cite[Sec.~3.2]{Drmota:1997:system-funct-equat}, we obtain that the
moment generating function is
\begin{equation*}
  \E(e^{\langle
    \bfOmega_n,\bfs\rangle})=\frac{y_{n}(e^{\bfs})}{y_{n}(\bfones)}
\end{equation*}
with $y_{n}(\boldsymbol{z})$ defined
in \cite{Drmota:1997:system-funct-equat}. Using
\cite[Equ.~(4.9)]{Drmota:1997:system-funct-equat}, this moment
generating function has an asymptotic expansion as in \eqref{eq:moment-asymp} with $\phi_{n}=n$. Thus
$\bfOmega_{n}$ is asymptotically normally distributed after
standardisation (as was shown in~\cite{Drmota:1997:system-funct-equat})
and additionally the speed of convergence is $O(n^{-1/2})$.

Other context-free languages can be analysed in the same way, either
by directly using the results in \cite{Drmota:1997:system-funct-equat} (if the
underlying system is strongly connected) or by similar methods. This
has applications, for example, in genetics (see~\cite{Poznanovic-Heitsch:2014:asymp-rna}).

\subsection{Dissections of Labelled Convex Polygons}
Let $S_{1}\cupdot\cdots\cupdot S_{t+1}=\{3,4,\ldots\}$ be a
partition. We dissect a labelled convex $n$-gon into smaller convex polygons by
choosing some non-intersecting diagonals. Each small polygon should be a
$k$-gon with $k\not\in S_{t+1}$. Define $a_{n}(\bfr)$ to be the number
of dissections of an $n$-gon such that it consists of exactly $r_{i}$ small
polygons whose number of vertices is in $S_{i}$, for $i=1$, \dots,
$t$. For convenience, we use $a_{2}(\bfr)=[\bfr=\bfzero]$. Asymptotic normality was proved in
\cite[Sec.~3]{Bender-Richmond:1983:centr}, see also
\cite[Ex.~7.1]{Bender:1974:asymp-method-enumer} for a one-dimensional
version. We additionally provide an estimate for the speed of convergence.

Let
\begin{equation*}
  f(z,\bfx)=\sum_{\substack{n\geq2\\ \bfr\geq 0}}a_{n}(\bfr)\bfx^{\bfr}z^{n-1}.
\end{equation*}
Then choosing a $k$-gon with $k\in S_{1}\cupdot\cdots\cupdot S_{t}$
and gluing dissected polygons to $k-1$ of its sides translates into
the equation
\begin{equation*}
  f=z+\sum_{i=1}^{t}x_{i}\sum_{k\in S_{i}}f^{k-1}.
\end{equation*}
Following \cite{Bender:1974:asymp-method-enumer}, this equation can be used to obtain an
asymptotic expression for the moment generating function as in
\eqref{eq:moment-asymp} with $\phi_{n}=n$. The asymptotic normal
distribution follows after suitable standardisation with speed of
convergence $O(n^{-1/2})$.

\section{Sketch of the Proofs}\label{sec:sketch-proofs}
We now sketch the main ideas of the proofs of
Theorems~\ref{theorem:Berry-Esseen-dimension-m}
and~\ref{th:quasi-power-dD}. All details can be found in the full
version of this extended abstract.

\begin{sketch}[Sketch of the proof of Theorem~\ref{theorem:Berry-Esseen-dimension-m}]
  As in
  \cite{Sadikova:1966:esseen:englisch,Gamkrelidze:1977:englisch,Rao:2002:anoth-esseen}, our
  proof of the Berry--Esseen inequality proceeds via adding a continuous random
  variable $\bfQ$ to
  our random variables $\bfX$ and $\bfY$. The
  characteristic function of $\bfQ$ vanishes outside $[-T, T]^m$. The error resulting from replacing the difference of the distribution functions $\abs{F_{\bfX}-F_{\bfY}}$ by $\abs{F_{\bfX+\bfQ} - F_{\bfY+\bfQ}}$ can be estimated by the final summand
  in~\eqref{eq:Berry-Esseen}. In principle, L\'evy's theorem then allows to
  bound the difference of the distribution functions by the difference of the 
  characteristic functions.  
   Instead of only using the difference of the
  characteristic functions, we use the difference $\abs{\Lambda(\varphi_{\bfX})-\Lambda(\varphi_{\bfY})}$, which
  ensures boundedness of the integral
  in~\eqref{eq:Berry-Esseen} at least if the first moments exist. However, we have to compensate $\Lambda$ by the
  sum over the differences of the marginal distribution functions,
  which yields the second summand in~\eqref{eq:Berry-Esseen}.
\end{sketch}

\begin{sketch}[Sketch of the proof of Theorem~\ref{th:quasi-power-dD}]
    First, the characteristic
   function of the standardised random variable $\bfX=(\bfOmega_{n}-\grad u(\bfzero)\phi_{n})/\sqrt{\phi_{n}}$ is
   \begin{equation*}
     \varphi_{\bfX}(\bfs)=\exp\Bigl(-\frac12 \bfs^\top \Sigma \bfs + O\Bigl(\frac{\norm{\bfs}^3+\norm{\bfs}}{\sphi}\Bigr)\Bigr)
   \end{equation*}
for $\norm{\bfs}<\tau\sqrt{\phi_{n}}/2$. Thus, we obtain convergence
in distribution as stated in the theorem. 

To obtain a bound for the speed of convergence, we use the Berry--Esseen
   inequality given in
   Theorem~\ref{theorem:Berry-Esseen-dimension-m} for $\bfY$ an
   $m$-dimensional normal distribution. We bound the
  difference of $\Lambda$ evaluated at the characteristic function of
  $\bfX$ and the one of the
  normal distribution by the exponentially decreasing function
\begin{equation*}
    \abs{\Lambda(\varphi_\bfX)(\bfs)-\Lambda(\varphi_\bfY)(\bfs)}\le
    \exp\Bigl(-\frac{\sigma}{4}\norm{\bfs}^2 + O(\norm{\bfs})\Bigr)O\Bigl(\frac{\norm{\bfs}^3+\norm{\bfs}}{\sphi}\Bigr)
  \end{equation*}
for suitable $\bfs$ where $\sigma$ is the smallest eigenvalue of $\Sigma$.

We then estimate the integral in
  \eqref{eq:Berry-Esseen}.
   For
  the variables in a neighbourhood of zero, we get rid of the denominator by
  Taylor expansion using the zero of $\Lambda(\varphi_{\bfX})-\Lambda(\varphi_{\bfY})$ at $\bfzero$. The error term of the Taylor expansion can be
  estimated by the difference of the characteristic
  functions  using Cauchy's formula. The
  exponentially small tails are used to bound the contribution of the
  large variables in the integral in~\eqref{eq:Berry-Esseen}.
  
  The second summand in
  \eqref{eq:Berry-Esseen} can be estimated inductively.
\end{sketch}

\bibliographystyle{amsplain}
\bibliography{bib/cheub}
\end{document}